\documentclass[a4paper,10pt]{amsart}
\usepackage{amsfonts}
\usepackage{enumerate}
\usepackage{amsmath, amsthm}
\usepackage{amscd}
\usepackage{amssymb}
\input xy
\xyoption{all}
\usepackage{latexsym,graphicx}
\usepackage{varioref}

\def\adots{\mathinner{\mkern2mu\raise 1pt\hbox{.}\mkern 3mu\raise
4pt\hbox{.}\mkern1mu\raise 7pt\hbox{{.}}}}

\newcommand{\sspace}{\vspace{0.25cm}}
\newcommand{\noi}{\noindent}

\theoremstyle{plain}
\newtheorem{theor}{Theorem}[subsection]

\newtheorem{prop}[theor]{Proposition}
\newtheorem{lem}[theor]{Lemma}

\newtheorem{defi}[theor]{Definition}

\theoremstyle{remark}
\newtheorem{rem}[theor]{Remark}

\newtheorem{Example}[theor]{\noi Example}

\numberwithin{equation}{section}

\newcommand{\CC}{{\mathbb C}}
\newcommand{\RR}{{\mathbb R}}

\newcommand{\ZZ}{{\mathbb Z}}

\newcommand{\G}{{\mathbf G}}
\newcommand{\HH}{{\mathbb H}}
\newcommand{\LL}{{\mathbf L}}

\newcommand{\Z}{{\mathbf Z}}

\newcommand{\NN}{{\mathbb N}}

\newcommand{\Ga}{\Gamma}

\newcommand{\lo}{\longrightarrow}

\newcommand{\End}{{\rm End}}

\newcommand{\Ad}{\textnormal{Ad} \;}

\newcommand{\GL}{{\rm \bf GL}}

\newcommand{\SU}{\mathbf{SU}}
\newcommand{\SO}{\mathbf{SO}}
\newcommand{\U}{\mathbf{U}}
\DeclareMathOperator{\HHom}{\textnormal{Hom}}

\newcommand{\M}{\mathbf{M}}

\newcommand{\Lieg}{\mathfrak{g}}
\newcommand{\LieZ}{\mathfrak{z}}

\newcommand{\Aut}{\mathbf{Aut}}
\newcommand{\VR}{V_{\RR}}
\newcommand{\VC}{V_{\CC}}

\newcommand{\VH}{V_{\HH}}
\newcommand{\Sp}{\mathbf{Sp}}
\newcommand{\LieSp}{\mathfrak{sp}}

\newcommand{\F}{\mathcal{F}}
\newcommand{\su}{\mathfrak{su}}
\newcommand{\LieU}{\mathfrak{u}}
\newcommand{\LieSO}{\mathfrak{so}}

\newcommand{\LieL}{\mathfrak{l}}
\newcommand{\LieK}{\mathfrak{k}}
\newcommand{\LieP}{\mathfrak{p}}

\newcommand{\HypH}{\mathbf{H}_\HH}
\newcommand{\HypC}{\mathbf{H}_\CC}
\newcommand{\hc}{h_\mathbb{C}}
\newcommand{\hh}{h_{\HH}}

\newcommand{\E}{\mathcal{E}}
\newcommand{\Lcal}{\mathcal{L}}

\newcommand{\gl}{\mathfrak{gl}}

\newcommand{\ww}{S^2 \VC^*}
\newcommand{\wwprime}{S^2 \LieP^*}

\newcommand{\bc}{{\mathbb C}}

\newcommand{\bh}{{\mathbb H}}

\newcommand{\br}{{\mathbb R}}

\newcommand{\ra}{\rightarrow}

\newcommand{\Gr}{\textnormal{Gr}}
\begin{document}
\title{Local quaternionic rigidity for complex hyperbolic lattices}
\author{I. Kim, B. Klingler and P. Pansu}

\date{}
\footnotetext[1]{I. Kim gratefully acknowledges the partial support
of KRF Grant (0409-20060066) and a warm support of IHES during his
stay.} \footnotetext[2]{B.Klingler gratefully acknowledges support of
NSF grant DMS 0635607} \footnotetext[3]{P. Pansu, Univ Paris-Sud, Laboratoire de Math\'ematiques d'Orsay, Orsay, F-91405 \hfil \break \indent CNRS, Orsay, F-91405.}

\maketitle

\begin{abstract}
Let $\Gamma \stackrel{i}{\hookrightarrow} L$ be a lattice in the real simple Lie group $L$.
If $L$ is of rank at least $2$ (respectively locally isomorphic to $Sp(n,1)$) any
unbounded morphism $\rho: \Gamma \lo G$ into a simple real Lie group
$G$ essentially extends to a Lie morphism $\rho_L: L \lo G$ (Margulis's superrigidity theorem, respectively
Corlette's theorem). In particular any such morphism is
infinitesimally, thus locally, rigid. 

On the other hand for $L=SU(n,1)$ even morphisms of the form $\rho :
\Gamma \stackrel{i}{\hookrightarrow} L \lo G$ are not infinitesimally
rigid in general. Almost nothing is known about their local rigidity. 
In this paper we prove that any {\em cocompact} lattice $\Ga$ in
$SU(n,1)$ is essentially locally rigid (while in general not infinitesimally
rigid) in the quaternionic groups $Sp(n,1)$, $SU(2n,2)$ or $SO(4n,4)$ (for the natural 
sequence of embeddings $SU(n,1) \subset Sp(n,1) \subset SU(2n,2)
\subset SO(4n,4))$.
\end{abstract}

\tableofcontents
\section{Introduction}

\subsection{Complex hyperbolic lattices and rigidity}
The main open question concerning lattices of Lie groups is certainly
the study of complex hyperbolic lattices and their finite dimensional
representations. Indeed, Margulis's super-rigidity theorem states that
any irreducible complex finite-dimensional representation of a lattice
$\Gamma$ of a simple real Lie group $L$ of real rank $r>1$ either has bounded image, or is the
restriction to $\Gamma$ of an irreducible finite-dimensional
representation of $L$. The remaining case of simple real
Lie groups of rank $1$ contains $3$ families~: the real hyperbolic
group $SO(n,1)$, the complex hyperbolic group $SU(n,1)$ and the
quaternionic hyperbolic group $Sp(n,1)$, plus one exceptional group
$F_4^{-20}$. Margulis's description has been extended to lattices of $Sp(n,1)$ and $F_4^{-20}$
by Corlette \cite{cor} and Gromov-Schoen \cite{grs}. On the other hand
one knows that $SO(n,1)$ admits
lattices with unbounded representations not coming from
$SO(n,1)$. Examples have been constructed by Makarov \cite{mak} and Vinberg
\cite{vin} for small $n$ and by Johnson-Millson \cite{jom} and Gromov- Piatetski-Shapiro
\cite{gps} for any $n \in \NN$. Concerning $SU(n,1)$, Mostow
\cite{mos} exhibited a striking counterexample
to superrigidity for $n=2$~:
namely two cocompact (arithmetic)
lattices $\Ga$ and $\Ga'$ in $SU(2,1)$ and a {\em surjective} morphism
$\rho: \Ga \lo \Ga'$ with {\em infinite kernel}.
Essentially nothing is known for $n >3$.

\sspace
In this paper, we restrict ourselves to the deformation theory of
complex hyperbolic cocompact lattices. Let $n>1$ be an integer and consider the complex hyperbolic group
$L=SU(n,1)$ : this is the group of real point of $\LL= \SU(n,1)=\SU (\VC, \hc)
$, the special unitary algebraic $\RR$-group of
linear isometries of $(\VC, \hc)$ where $\VC$ denotes the $(n+1)$-dimensional
$\CC$-vector space endowed with the Hermitian form
$\hc(\mathbf{z},\mathbf{w}) =-z_0 \overline{w_0} + z_1 \overline{w_1}
\cdots  + z_n \overline{w_n}$. Let $i: \Gamma \hookrightarrow SU(n,1)$
be a cocompact complex hyperbolic lattice. Let $j : \SU(n,1) \hookrightarrow \G$ be an injective
$\RR$-morphism of $\RR$-algebraic groups. Does there exist any {\em non-trivial} deformation
of $\rho = j \circ i : \Gamma \lo G = \G(\RR)$, i.e. a
continuous family of morphisms $\rho_t : \Ga \lo G$, $t \in I = [0,1]$, with $\rho_0 =
\rho$ not of the form $\rho_t = g_t \cdot \rho \cdot g_t^{-1}$ for
some continuous family $g_t \in G$, $t \in I$~?

\subsection{First order deformations}
Let $\M(\Ga, \G)(\RR)= (\bf{\textnormal{Hom}} (\Ga, \G)//\G) (\RR)$ be the moduli space of representations of $\Ga$
in $\G(\RR)$ up to conjugacy.
The space of first-order deformations of $\rho$,
i.e. the real Zariski tangent space at $[\rho]$ to $\M(\Ga, \G)(\RR)$,
naturally identifies with the first cohomology group $H^1(\Ga, \Ad
\rho)$, where $\Ad \rho : \Ga \stackrel{\rho}{\hookrightarrow} G \stackrel{\Ad}{\rightarrow}
\text{Aut}(\Lieg)$ is the natural representation deduced from $\rho$
and the adjoint action of $G$ on its Lie algebra $\Lieg$. Thus the
non-vanishing of $H^1(\Ga, \Ad \rho)$ is a necessary condition for
$\M(\Ga, \G)(\RR)$ not being trivial at the point $[\rho]$. Raghunathan
\cite{rag} gave the list of irreducible finite-dimensional
$\SU(n,1)$-modules which may have non-vanishing $\Ga$-cohomology in
degree~$1$~:

\begin{theor}[Raghunathan]
Let $\lambda : \SU(n,1) \lo \GL(W)$ be a real
finite dimensional irreducible representation of $\SU(n,1)= \SU(\VC,\hc)$. Let $\Ga$
be a cocompact lattice in $SU(n,1)$. Then $H^1(\Ga, W)=0$ except if
$W\simeq  S^j\VC$ for some $j \geq 0$, where $S^j$
denotes the $j-$th symmetric power.
\end{theor}

\begin{rem}
In this theorem $\VC$ is seen as a {\em real} representation. In
particular $S^j \VC^* \simeq S^j\VC$ as a real $SU(n,1)$-module.
\end{rem}

As a corollary, $[\rho]\in \M(\Ga, \G)(\RR)$ is isolated except maybe if
$\Ad j: \SU(n,1) \lo \Aut(\Lieg)$ contains an
$\SU(n,1)$-direct factor isomorphic to $S^j\VC$ or $S^j \VC^*$ for some
integer $j \geq 0$.

\rem
For each $n$ and each $j$ one can, following a method first introduced
by Kazhdan,
exhibit a cocompact lattice $\Ga$ of $SU(n,1)$ such that $H^1(\Ga,
S^j\VC) \not = 0$, c.f. \cite[chap. VIII]{bw}.

\begin{Example}
Let $\Ga \stackrel{i}{\hookrightarrow} SU(n,1)$ be a cocompact
  lattice. Let $j = \text{Id} : \SU(n,1) \lo \SU(n,1)$. By Raghunathan's theorem,
  $H^1(\Ga, \Ad i)=0$, thus $\Ga$ cannot be non-trivially deformed in
  $SU(n,1)$. This was already proved by Weil \cite{wei}.
\end{Example}

\begin{Example} \label{ortho}
Let $j : \SU(n,1)= \SU(\VC,\hc) \hookrightarrow \mathbf{SO}(2n,2)=
\mathbf{SO}((\VC)^\RR, \text{Re} \, \hc)$ be the natural
embedding. Notice that $j$ factorizes as $\SU(n,1) \hookrightarrow \U(n,1)
\hookrightarrow \mathbf{SO}(2n,2)$. One easily checks that the Lie
algebra $\mathfrak{so}(2n,2)$ is isomorphic
as an $\SU(n,1)$-module to the direct sum of irreducible modules $\RR
\oplus \mathfrak{su}(n,1) \oplus \Lambda^2 \VC$, where $\RR =
\text{Lie}(\Z(\RR))$ is the Lie algebra of the centralizer $\Z$ of
$\SU(n,1)$ in $\U(n,1)$.
Thus $H^1(\Ga, \Ad \rho) = H^1(\Ga, \RR)$ and any deformation of $\rho$ in
$SO(2n,2)$ is of the form $\rho \cdot \chi$, where
$\chi : \Ga \lo \Z(\RR) = S^1$ is a unitary character of $\Ga$.
\end{Example}

\subsection{Local rigidity}

\subsubsection{Formal completion of $\M(\Ga, \G)(\RR)$ at $[\rho]$}
Studying first-order deformations is not enough the local rigidity
problem stated in the introduction~:
even if $j : \SU(n,1) \hookrightarrow \G$ is such that {a priori} $H^1(\Ga,
\Ad \rho)$ does not vanish it may happen that very few of these infinitesimal
deformations can be integrated. However it is enough to study second
order deformations. Let \- $\HypC^n = SU(n,1)/U(n)$ denote
the symmetric space of $SU(n,1)$~: this is the complex hyperbolic $n$-space of
negative lines in $(\VC, \hc)$), it is naturally endowed with an $SU(n,1)$-invariant
K\"ahler form $\omega_{\HypC^n}$. Without loss of generality (passing
to a finite index subgroup) one can assume that $\Ga$ is torsion-free,
so that $M= \Ga \backslash \HypC^n$ is a compact K\"ahler manifold with
fundamental group $\Ga$. One can then apply the following formality theorem of Goldman-Millson
\cite{gom} (for the case of complex variations of Hodge structures)
and Simpson \cite{sim1} (in general)~:

\begin{theor}[Goldman-Millson, Simpson]
Let $M$ be a connected compact K\"ahler manifold with fundamental group
$\Ga$, $\G$ a real reductive algebraic group and $\rho: \Ga \lo G=\G(\RR)$ a reductive representation. Let $C \subset
H^1(\Ga, \Ad \rho)$ be the affine cone defined by
$$ C = \{ u \in H^1(\Ga, \Ad \rho) \; / [u, u] = 0 \in H^2(\Ga, \Ad
\rho) \} \;\;.$$
Then the formal completion of $\M(\Ga, \G)(\RR)$ at $[\rho]$ is
isomorphic to the formal completion of the good quotient $C/H$, where
$H$ denotes the centralizer of $\rho(\Ga)$ in $G$.
\end{theor}

\subsubsection{Goldman-Millson rigidity result}
The first result about non-integra\-bility of some first-order
deformations for cocompact complex hyperbolic lattices
is due to Goldman-Millson \cite{gm}~: they consider the embedding
$$j : \SU(n,1)= \SU(\VC, \hc) \hookrightarrow \SU(n+1,1)= \SU(\VC \oplus \CC,
\hc \oplus 1)\;\;.$$
In this case the space of first-order deformations $H^1(\Ga, \Ad
\rho)$ at $\rho = j \circ i$ decomposes as $H^1(\Ga, \RR) \oplus
H^1(\Ga, \VC)$. The first summand $H^1(\Ga, \RR)$
corresponds once more to the uninteresting deformations obtained by deforming
$\Ga$ in $U(n,1)$ by a curve of homomorphism into the centralizer $\Z= \U(1)$ of
$\SU(n,1)$ in $\U(n,1)$. The second summand, which potentially corresponds to
Zariski-dense deformations of $\rho$ in $SU(n+1,1)$, is non-zero for
general $\Gamma$. However Goldman and Millson prove that none of these
deformations can be integrated. Thus any representation $\lambda: \Ga
\lo SU(n+1,1)$ sufficiently close to $\rho$ is conjugate to a
representation of the form $\rho
\cdot \chi$, where $\chi : \Ga \lo Z=S^1$.
A similar result can be obtained by replacing the natural embedding $j :
\SU(n,1) \hookrightarrow \SU(n+1,1)$ with the natural embedding $j :
\SU(n,1) \hookrightarrow \SU(n+k,1)$ for some integer $k \geq 1$.

\subsubsection{Possible extensions}
One natural generalization of Goldman-Millson's result consists in
studying {\em global rigidity} of representations $\rho: \Ga \lo G$ with $G$ simple of Hermitian
type, under certain assumptions on $\rho$. Let $X_G$ be the (K\"ahler) symmetric space associated to $G$,
with K\"ahler form $\omega_G$. Let $\omega_M$ be the natural K\"ahler
form on $M$. Let $f: \tilde{M}= \HypC^n \lo X_G$ be any smooth $\rho$-equivariant
map. The de Rham class $[f^*\omega_G] \in H^2_{dR}(M)$ depends
only on $\rho$, not on $f$, and will be denoted $[\rho^* \omega_G]$.
Define the {\em Toledo invariant} $\tau(\rho)$ of $\rho$ as the
number
$$ \tau(\rho) = \frac{1}{n !} \int_M \rho^*\omega_G \wedge
\omega_M^{n-1} \;\;.$$
One easily shows that $\tau$ is a locally constant function on
$\M(\Ga, \G)(\RR)$. Moreover it satisfies a Milnor-Wood inequality~: under
suitable normalizations of the metrics one has
$$ |\tau(\rho)| \leq \text{rk} \, X_G \cdot \text{Vol}(M) \;\;.$$
One expects a {\em global rigidity result} for representations $\rho :
\Ga \lo G$ with {\em maximal Toledo invariant}~: namely $\rho$ is
expected to be faithful, discrete and stabilizing a holomorphic totally geodesic
copy of $\HypC^n$ in $X_G$. This has been proven by Corlette
\cite[theor. 6.1]{cor0} when $G$ is of rank
one and $\Ga$ cocompact (thus
generalizing Goldman-Millson's result), then by
B\"urger-Iozzi \cite{bi} and Koziarz-Maubon \cite{km} for
$G$ of rank $1$ and any complex hyperbolic lattice $\Ga$. Recently
Koziarz-Maubon \cite{km2} proved it when the group $G$ is of real rank
$2$. In the same kind of direction, we also refer to \cite{cor1/2}.

\subsection{The main result}
From the point of view of non-abelian Hodge theory, it is natural to
enlarge the study of representations of complex hyperbolic lattices
into groups of Hermitian type to the study of representations into groups
of Hodge type (i.e. simple real Lie groups admitting discrete series).
Among groups of Hodge type there is a particularly simple subclass~: the groups of {\em quaternionic
  type}, that is such that the associated symmetric space $X_G$ is quaternionic-K\"ahler.
The classical families in this class are $Sp(n,1)$, $SU(n,2)$ and
$SO(n,4)$, $n \geq 1$. The corresponding $3$ families of quaternionic K\"ahler non-compact
irreducible symmetric spaces of dimension $4n$, $n \geq 2$, are: $\HypH^n = Sp(n,1)/Sp(n) \cdot Sp(1)$,
$X^n= SU(n,2)/S(U(n) \times U(2))$ and $Y^n =SO(n,4)/S(O(n) \times O(4))$.
The only K\"ahler ones are $X^n$ and $Y^2$.

\sspace
The main result of this paper study quaternionic deformations of
cocompact complex hyperbolic lattices.
Let $\VH = \VC \otimes_{\CC} \HH $ be the
quaternionic right vector space of dimension $n+1$ (thus of real
dimension $4n+4$) endowed with the quaternionic Hermitian form
$\hh$ of signature $(n,1)$ deduced from
$\hc$. The complex Hermitian part $H$ of $\hh$ is a complex Hermitian
form on $\VC \oplus j \VC$ of signature $(2n,2)$.
Let $\Sp(n,1) = \SU(\VH, \hh)$ be the special unitary algebraic $\RR$-group of
linear transformation of $(\VH, \hh)$, $\U(2n,2)$ the unitary
$\RR$-group of linear transformations of $(\VC \oplus j \VC, H)$ and
$\SO(4n,4)$ the special orthogonal group of linear transformation of
$((\VH)_{\RR}, \textnormal{Re} H)$.
One obtains a natural sequence of embeddings
$$ \SU(n,1) \stackrel{j_{\U(n,1)}}{\hookrightarrow} \U(n,1)
\stackrel{j_{\Sp(n,1)}}{\hookrightarrow} \Sp(n,1)
\stackrel{j_{\U(2n,2)}}{\hookrightarrow} \U(2n,2)
\stackrel{j_{\SO(4n,4)}}{\hookrightarrow} \SO(4n,4)$$
corresponding to equivariant totally geodesic embeddings of symmetric spaces
$$ \HypC^n \stackrel{f_{\HypH^n}}{\hookrightarrow} \HypH^n
\stackrel{f_{X^{2n}}}{\hookrightarrow} X^{2n} \stackrel{f_{Y^{4n}}}{\hookrightarrow}
Y^{4n} \;\;.$$

\begin{rem} \label{rem1}
Notice that the totally geodesic embedding $\HypC^n \stackrel{f_{X^{2n}}
  \circ f_{\HypH^n}}{\hookrightarrow} X^{2n}$ between Hermitian symmetric spaces {\em is not
  holomorphic}~: the pull-back $(f_{X^{2n}} \circ f_{\HypH^n})^* \omega_{X^{2n}}$ is
identically zero.
\end{rem}

\sspace
For $i: \Ga \hookrightarrow SU(n,1)$ a cocompact lattice, and $\G= \U(n,1)$, $\Sp(n,1)$,
$\U(2n,2)$ or $\SO(4n,4)$  let $\rho_\G : \Ga \lo G$ be the composition $j_\G \circ \cdots \circ j_{\U(n,1)} \circ
i$. The space of first-order deformations $H^1(\Ga, \Ad
\rho_\G)$ at $[\rho_\G]$ is non-trivial for general $\Ga$. As in
Goldman-Millson's result we however prove~:

\begin{theor} \label{main_theorem}
Let $\Ga \stackrel{i}{\hookrightarrow} SU(n,1)$ be a cocompact
lattice and $\G$ one of the groups $\Sp(n,1)$,
$\U(2n,2)$ or $\SO(4n,4)$. Then any morphism $\lambda : \Ga \lo G=\G(\RR)$ close enough
to $\rho_\G$ is conjugate to a representation of the form $\rho_\G \cdot
\chi$, where $\chi : \Ga \lo Z_G(SU(n,1))$ (thus $Z_{Sp(n,1)}(SU(n,1)) =U(1)$ and
$Z_{U(2n,2)}(SU(n,1)) = Z_{SO(4n,4)}(SU(n,1)) = U(1) \times U(1)$).
\end{theor}

\begin{rem} \label{rem2}
Following remark~\ref{rem1} notice that the representation
$\rho_{\U(2n,2)}: \Ga \lo U(2n,2)$ satisfies
$\tau(\rho_{\U(2n,2)})=0$, thus has the smallest possible (in
absolute value) Toledo invariant. In particular
theorem~\ref{main_theorem} in this case is not covered by
Koziarz-Maubon \cite{km2} (nor Corlette \cite{cor1/2}). Also the
same method applies to prove the case when $G= U(n+k,m)$, more
generally when $G=Sp(n+k,m)$.
\end{rem}

\subsection{Organization of the paper}
The proof of theorem~\ref{main_theorem} essentially reduces to the case $\G
= \Sp(n,1)$, with an extra argument for $\SO(4n,4)$ (c.f. section~\ref{infi}).
In sections~\ref{classic}, \ref{conseq} and ~\ref{second}, we give a
first proof of the main theorem~\ref{main_theorem} using Goldman-Millson's strategy~:  first, using Matsushima and
Murakami's method \cite{matmu}, we show that harmonic 1-forms
representing nontrivial classes in $H^{1}(\Gamma,\Ad \rho_{\Sp(n,1)})$
are severely restricted : most of their components vanish, and one can
interpret them as $(1,0)$-forms $\alpha$ with values in a certain
complex vector bundle. Then we show that the cup-square
$[\alpha, \alpha] \in H^2(\Ga,\Ad \rho_{\Sp(n,1)}) $ paired with the K\"ahler form of complex
hyperbolic space is proportional to the squared $L^2$-norm of
$\alpha$, which implies the result.

In sections~\ref{nonab}, we indicate a more geometric proof of the
main theorem~\ref{main_theorem} based on
period domains and a result of Carlson-Toledo \cite{cto}.

\section{Infinitesimal deformations of lattices of $SU(n,1)$ in $G$}
\label{infi}

\subsection{The groups}

\begin{defi} Let $n>1$ be an integer.
We denote by $\VR$ the $n+1$-dimensional $\RR$-vector space, $\VC =
\VR \otimes_\RR \CC$ its complexification and $\VH = \VR\otimes_{\RR}
\HH $ its quaternionification (thus $\VH$ is a
{\em right} quaternionic vector space). We define $\GL(n+1, \HH)$ as
the $\RR$-group of $\HH$-linear automorphism of $\VH$.
\end{defi}

\begin{defi}
Let $Q_{\RR}$ be a real quadratic form of signature $(n,1)$
on $\VR$. We denote by $Q_\CC$ (respectively $Q_{\HH}$) its
complexification (resp. its quaternionification) on $\VC$ (resp. on
$\VH$).
\end{defi}

\begin{defi}
We denote by $\hc$ the complex Hermitianisation of $Q_\RR$ on $\VC$
and by $\hh$ the quaternionic Hermitianisation of $Q_\RR$ on $\VH$.
Thus $\hc(\mathbf{z}, \mathbf{w})= Q_{\CC}(\mathbf{z}, \overline{\mathbf{w}}^\CC)$ where
$\overline{\mathbf{w}}^\CC$ denotes the complex conjugate of $w \in \VC$ and
$\hh(\mathbf{z}, \mathbf{w})= Q_{\HH}(\mathbf{z}, \overline{\mathbf{w}}^\HH)$ where
$\overline{\mathbf{w}}^\HH$ denotes the quaternionic conjugate of $w \in \VH$.
\end{defi}

\noi
On the complex vector space $\VH = \VC \oplus j\VC$, the quaternionic Hermitian
form $h_{\HH}(\mathbf{z}, \mathbf{w})$ can be written as
$$h_\HH(\mathbf{z}, \mathbf{w}) = H(\mathbf{z}, \mathbf{w}) - j
\Omega(\mathbf{z}, \mathbf{ w})\;\;,$$
where $H$ is a complex Hermitian form on $\VC \oplus j \VC$ and $\Omega$ is the skew-symmetric
complex bilinear form on $\VC \oplus j \VC$ defined by
$\Omega(\mathbf{z}, \mathbf{w})=  H(\mathbf{z} \cdot j, \mathbf{\bar w})$.

\begin{defi}
We define the real algebraic groups~:
\begin{itemize}
\item  $\Sp(n,1)= \Sp(\VH, \hh)$ as the subgroup of $\GL(n+1, \HH)$
preserving $\hh$.
\item $\U(2n,2)$ the unitary group $\U(\VC \oplus
j\VC, H)$.
\item $\SO(4n,4)$ the special orthogonal group $\SO((\VC \oplus
j\VC)_\RR, \textnormal{Re} H)$.
\end{itemize}
Moreover we denote by $\Sp(2n+2, \CC)$ the complex symplectic group $\Sp(\VC \oplus j \VC, \Omega)$.
\end{defi}

The previous discussion implies immediately (where we consider
$\Sp(2n+2, \CC)$ as a real algebraic group)~:
\begin{lem}
$\Sp(n,1) =  \GL(n+1,\HH) \cap \U(2n,2) = \Sp(2n+2, \CC) \cap \U(2n,2)$.
\end{lem}

Consider the sequence of natural embeddings~:
\begin{equation} \label{sequence}
\SU(n,1) \stackrel{j_{\U(n,1)}}{\hookrightarrow} \U(n,1)
\stackrel{j_{\Sp(n,1)}}{\hookrightarrow} \Sp(n,1)
\stackrel{j_{\U(2n,2)}}{\hookrightarrow} \U(2n,2)
\stackrel{j_{\SO(4n,4)}}{\hookrightarrow} \SO(4n,4)\;\;.
\end{equation}

\begin{lem}
The sequence~(\ref{sequence}) induces an exact sequence of
$U(n,1)$-modules (under the adjoint representation)~:
\begin{equation} \label{exact}
\begin{split}
0 \lo \LieU(n,1) \lo &\LieSp(n,1) = \LieU(n,1) \oplus S^2 \VC^* \lo
\LieU(2n,2) = 2 \LieU(n,1) \oplus \Lambda^2 V_{\bc}^*  \oplus S^2 \VC^*
\lo \\
&\lo \LieSO(4n,4)= 2 \LieU(n,1) \oplus 2\Lambda^2 \VC^* \oplus 2\Lambda^2
\bar{V}_{\bc}^* \oplus S^2 \VC^* \oplus S^2 \bar{V}_{\bc}^*
\;\;.
\end{split}
\end{equation}
\end{lem}

\begin{proof}

{\em Case $\G=\Sp(n,1)$}. Let $M\in \mathfrak{gl}(n+1,\bh)$ and $M=C+jD$ where $C,D\in
\mathfrak{gl}(n+1,\bc)$. Then $M\in \mathfrak{sp}(n,1)$, if and only
if $C\in \mathfrak{u}(n,1)$ and $JD$ is symmetric where $J$ is the
diagonal matrix with entries $1,\cdots, 1,-1$. Write $E=JD$. If $A\in U(n,1)$,
\begin{eqnarray*}
A^{-1}MA=A^{-1}(C+jJE)A=A^{-1}CA+jJA^{\bot}EA.
\end{eqnarray*}
So, under $U(n,1)$,
$\mathfrak{sp}(n,1)=\mathfrak{u}(n,1)\oplus S^2 \VC^{*}.$

\sspace
{\em Case $\G=\SU(2n,2)$}. Let $q=a+jb$ be a quaternion, with $a$,
$b\in\bc$. The matrix of left multiplication by $q$ is
$\left (\begin{smallmatrix}a&-\bar{b}\\ b&\bar{a}\end{smallmatrix} \right )$. Therefore, if
$A\subset GL(n+1,\bc)$, its image under the embeddings $GL(n+1,\bc)\to
GL(n+1,\bh)\to GL(2n+2,\bc)$ is $\left (  \begin{smallmatrix}A&0\\
  0&\bar{A}\end{smallmatrix} \right )$. If
$\mathfrak{u}(n,1)\subset\mathfrak{gl}(n+1,\bc)$ is the subspace of
matrices $A$ such that $A^* Q+QA=0$ with $Q=\left (  \begin{smallmatrix}I_n&0\\
  0&-1\end{smallmatrix} \right )$, then $\LieU(n,1)$ is mapped to
$\mathfrak{u}(2n,2)$ defined as the subspace of matrices
$M\in\mathfrak{gl}(2n+2)$ such that $M^* Q'+Q'M=0$ with
$Q'=\left (  \begin{smallmatrix}Q&0\\ 0&Q\end{smallmatrix} \right )$. Thus, under the adjoint
action of $U(n,1)$,
\begin{eqnarray*}
\LieU(2n,2)=\LieU(n,1)\oplus \LieU(n,1)\oplus Hom_{\bc}(\bc^{n+1},\bc^{n+1}),
\end{eqnarray*}
where $U(n,1)$ acts on a square matrix $N\in Hom_{\bc}(\bc^{n+1},\bc^{n+1})$ as follows,
\begin{eqnarray*}
(A,N)\mapsto A^{-1}N\bar{A}.
\end{eqnarray*}
Putting $B=NQ$ conjugates this action to
\begin{eqnarray*}
(A,B)\mapsto A^{-1}B(A^{-1})^{\top},
\end{eqnarray*}
i.e. $Hom_{\bc}(\bc^{n+1},\bc^{n+1})=V_{\bc}^* \otimes_{\bc} V_{\bc}^*$. Then
\begin{eqnarray*}
\LieU(2n,2)=2\LieU(n,1)\oplus\Lambda^2 V_{\bc}^* \oplus S^2 V_{\bc}^* ,
\end{eqnarray*}
where $S^2 V_{\bc}^* $ corresponds to matrices of the form
$\left (  \begin{smallmatrix}0&BQ\\ -B^* Q&0\end{smallmatrix} \right )$ in $\LieU(2n,2)$ with $B$
symmetric.

\sspace
{\em
Case $\G=\SO(4n,4)$}. We have seen that the embedding $GL(n+1,\bc)\to
GL(n+1,\bh)\to GL(2n+2,\bc)$ lands into the block diagonal subgroup
$GL(n+1,\bc)\times GL(n+1,\bc)\subset GL(4n+4,\br)$. In particular,
$U(n,1)$ lands into $U(n,1)\times U(n,1)\subset O(2n,2)\times
O(2n,2)$. Under $O(2n,2)\times O(2n,2)$,
\begin{eqnarray*}
\mathfrak{so}(4n,4)=\mathfrak{so}(2n,2)\oplus \mathfrak{so}(2n,2)\oplus End_{\br}(\br^{2n+2}).
\end{eqnarray*}
Since $U(n,1)$ preserves a complex structure, $\br^{2n+2}=\bc^{n+1}$,
every $\br$-linear map $L$ is the sum of a $\bc$-linear and an
anti-$\bc$-linear one, $L=L_{\bc}+L_{\bar{\bc}}$, and the action of
$A\in U(n,1)$ on $L$ is
$A^{-1}L_{\bc}\bar{A}+\bar{A}^{-1}L_{\bar{\bc}}A$. Thus
$End_{\br}(\bc^{n+1})$ equals the sum of
$End_{\bc}(\bc^{n+1})=V_{\bc}^* \otimes_{\bc} V_{\bc}^*$ and its
conjugate $\bar{V}_{\bc}^{*}\otimes_{\bc}\bar{V}_{\bc}^{*}$.

The map $\mathfrak{so}(2n,2)\to \Lambda^2 (\br^{2n+2})^*$, $C\mapsto
QC$ conjugates the adjoint $SO(2n,2)$ action with its action on real
alternating 2-forms. In presence of the $U(n,1)$-invariant complex
structure $J$, alternating 2-forms split into two subspaces
$\Lambda_+$ and $\Lambda_-$. Indeed, $\Lambda^2 J$ is an
involution. The inverse map $B\mapsto QB$ maps $\Lambda_+$ to
$\mathfrak{u}(n,1)\subset\mathfrak{so}(2n,2)$. $J$ also acts as a
derivation on alternating 2-forms, yielding a complex structure on
$\Lambda_-$. Since
\begin{eqnarray*}
\mathfrak{so}(2n,2)\otimes\bc =\Lambda^{2,0}(\br^{2n+2})^* \oplus\Lambda^{1,1}(\br^{2n+2})^* \oplus\Lambda^{0,2}(\br^{2n+2})^* ,
\end{eqnarray*}
$\Lambda^{1,1}(\br^{2n+2})^* =\Lambda_+ \otimes\bc$,
$\Lambda^{2,0}(\br^{2n+2})^* \oplus\Lambda^{0,2}(\br^{2n+2})^*
=\Lambda_-\otimes\bc$, thus, as a complex representation of $U(n,1)$,
the $\Lambda_-$ factor in the first diagonal block is isomorphic to
$\Lambda^{2}V_{\bc}^*$, and the $\Lambda_-$ factor in the second
diagonal block is isomorphic to $\Lambda^{2}\bar{V}_{\bc}^{*}$.

We conclude that
\begin{eqnarray*}
\mathfrak{so}(4n,4)=\mathfrak{z}\oplus
2\mathfrak{su}(n,1)\oplus 2\Lambda^{2}V_{\bc}^* \oplus 2\Lambda^{2}\bar{V}_{\bc}^{*}\oplus S^2 V_{\bc}^* \oplus S^{2}\bar{V}_{\bc}^{*},
\end{eqnarray*}
where $\mathfrak{z}=\br^2$ is the sum of the centers of the 2 copies of $\mathfrak{u}(n,1)$, generated respectively by $\left(\begin{smallmatrix}J&0\\ 0&0\end{smallmatrix}\right)$ and $\left(\begin{smallmatrix}0&0\\ 0&J\end{smallmatrix}\right)$.
\end{proof}

Choose $J=\left(\begin{smallmatrix}
0 & -I_{n+1}\\
I_{n+1} & 0
\end{smallmatrix}\right)$ as a complex structure on $\br^{2n+2}$.

\begin{lem}If $M=A+iB$ is a complex matrix representing an
anti-$\bc$-linear map, it is mapped to $\left(\begin{smallmatrix}
                                          A & B \\
                                          B & -A
                                          \end{smallmatrix}\right)$
in $GL(2n+2, \br)$.

Then $Z\in S^2 V_{\bc}^* , Z'\in S^2 \bar{V}_{\bc}^*$ can be
written
$$
Z=\begin{pmatrix}0&BQ'\\ -B^*Q'&0\end{pmatrix},\quad Z'=\begin{pmatrix}0&B'Q'\\ -B'^*
Q'&0\end{pmatrix}
$$
respectively, where $B=\left(\begin{smallmatrix} C & -D \\ D & C
\end{smallmatrix}\right),
 B^*=\left(\begin{smallmatrix} C & D \\ -D & C\end{smallmatrix}\right),
B'=B'^* =\left(\begin{smallmatrix} C' & D' \\ D' & -C'\end{smallmatrix}\right)$ and
$C$, $D$, $C'$, $D'$ are symmetric real matrices.
\end{lem}

\begin{proof}
The first statement comes directly from calculation.
The second follows from the first statement and the fact that the
matrices are symmetric and the fact that they are in
$\mathfrak{so}(4n,4)$.
\end{proof}

\subsection{Some reductions}

\begin{lem}
The special case of the main theorem~\ref{main_theorem} for $\G = \Sp(n,1)$
implies the main theorem for $\G = \SU(2n,2)$, but not quite for $\G = \SO(4n,4)$ .
\end{lem}

\begin{proof}
One deduces from the sequence~(\ref{sequence}) the following commutative
diagram~:
$$
\xymatrix{
H^1(\Ga, \Ad \rho_{\Sp(n,1)}) \ar[d]_q \ar@{^(->}[r]^{j_{\U(2n,2)}} &
H^1(\Ga, \Ad \rho_{\U(2n,2)}) \ar[d]_q \ar@{^(->}[r]^{j_{\SO(4n,4)}} &
H^1(\Ga, \Ad \rho_{\SO(4n,4)}) \ar[d]_q \\
H^2(\Ga, \Ad \rho_{\Sp(n,1)}) \ar@{^(->}[r]_{j_{\U(2n,2)}} &
H^2(\Ga, \Ad \rho_{\U(2n,2)}) \ar@{^(->}[r]_{j_{\SO(4n,4)}} &
H^2(\Ga, \Ad \rho_{\SO(4n,4)})
}
\;\;,$$
where $q: H^1(\Ga, \Ad \rho_G) \lo H^2(\Ga, \Ad \rho_G)$ denotes the
quadratic map deduced from the symmetric
bilinear map $$[\cdot, \cdot]: H^1(\Ga, \Ad \rho_G) \times H^1(\Ga, \Ad
\rho_G) \lo H^2(\Ga, \Ad \rho_G)\;\;.$$
As
$H^1(\Ga, \Lambda^2 \VC^*)= H^1(\Ga, \su(n,1))=0$ and
as the space $H^1(\Ga, z_{\Lieg}(\su(n,1))$ belongs to the null-space
of the quadratic map $q$, the proof of the main theorem for $\G =\Sp(n,1)$
or $\SU(2n,2)$ reduces to showing that the quadratic map
$q: H^1(\Ga, S^2 \VC^*) \lo H^2(\Ga, \LieSp(n,1)) \subset H^2(\Ga,
\LieU(2n,2))$ is anisotropic. Thus solving the case $\G =
\Sp(n,1)$ simultaneously solves the case $\G = \SU(2n,2)$. However, the
proof of the main theorem for $\G =\SO(4n,4)$, which amounts to showing that the
quadratic map
$$q: H^1(\Ga, S^2 \VC^*) \oplus H^1(\Ga, S^2 \bar{V}_{\bc}^*)\lo H^2(\Ga, \LieSO(4n,4))$$
is anisotropic, requires an extra computation.
\end{proof}

\section{A classical vanishing theorem} \label{classic}

\subsection{Matsushima and Murakami's vanishing theorem}
\label{MM}

Let $\LL$ be a simple real algebraic group of non-compact type, $L =
\LL(\RR)$ its Lie group of real points, $K$ a maximal compact subgroup
of $L$, $\theta: \LieL \lo \LieL$ the Cartan involution associated to $K$ of the Lie
algebra $\LieL$ of $L$, $\LieL = \LieK \oplus \LieP$ the
Cartan decomposition associated to $\theta$, $X = L/K$ the symmetric
space of $L$.

Let $i: \Gamma \hookrightarrow L$ be a (torsion free) cocompact lattice and $p:
\Gamma \backslash L \lo M = \Gamma \backslash L /K$ the natural
principal $K$-bundle on the locally symmetric manifold $M$.
Let $\rho : \LL \lo \GL(F)$ be a finite dimensional representation of
$\LL$. For $p$ a positive integer, the cohomology $H^p(\Ga, F)$ is
canonically isomorphic to the cohomology $H^p(M, F_\rho)$ of the local
system $F_\rho$ on $M$ associated to $\rho$, which can be computed
using the usual de Rham complex $(C^\bullet(M, F_\rho), d)$.

\sspace
Fix an admissible inner product $(,)_F$ on $F$, i.e. one which is
$\rho(K)$-invariant and for which elements of $\rho(\mathfrak{p})$ are
symmetric. This is enough to define a natural Laplacian $\Delta:
C^\bullet(M, F_\rho) \lo C^\bullet (M, F_\rho)$ and prove that $H^p(M,
F_\rho)$ is isomorphic to the space $$\mathcal{H}^p(M, F_{\rho})= \{\eta
\in C^P(M, F_\rho)\;/\; \Delta \eta =0\}$$ of harmonic forms
\cite[section 6]{matmu}.

\sspace
Following p. 376 of \cite{matmu}, define an $F$-valued differential form $\eta^0$ on $G$ as follows.
\begin{eqnarray}
\eta^0_s=\rho(s^{-1})\pi^* \eta_s, \quad s\in G.\end{eqnarray}
Fix a Killing-orthonormal basis $X_1,\ldots,X_N$ of
$\mathfrak{p}$. The induced inner product on $\HHom(\mathfrak{p},F)$ is
given by
$$(\eta,\zeta)=\sum_{h=1}^N (\eta(X_h),\zeta(X_h))_F .$$

\begin{defi}
Let $p$ be a positive integer. One defines a symmetric operator
$T_p$ on $\HHom(\mathfrak{p},F)$ as follows.
$$
\forall \, \eta \in \HHom(\mathfrak{p},F), \;\; \forall \, Y \in \LieP, \;\;
T_p\eta(Y)=\frac{1}{p} \sum_{k=1}^N \rho(X_k)^2 \eta(Y)+\rho([Y,X_k])\eta(X_k)\;\;.
$$
\end{defi}

\begin{theor}[Matsushima-Murakami] \cite[theor.7.1]{matmu}
\label{matu}
If $\eta$ is a harmonic $p$-form on $M=\Gamma\setminus G/K$, then
$$\int_{\Gamma\setminus G}(T\eta^0 ,\eta^0 )\leq 0.$$
As a consequence, if the symmetric operator $T_p$ on $\HHom(\mathfrak{p},F)$ is positive definite,
then the cohomology
group $H^p (\Gamma,F_\rho)$ vanishes.
\end{theor}

\subsection{Case of 1-forms}

\begin{prop}
\label{calculT}
Let $\eta\in\HHom(\mathfrak{p},F)$. Let $\beta: \LieP \otimes \LieP \lo
F$ denote the $F$-valued
bilinear form on $\mathfrak{p}$ defined by
$
\beta(X,Y)=\rho(X)(\eta(Y)).
$
Split $\beta=\sigma+\alpha$ into its symmetric and skew-symmetric parts. Then
$
(T\eta,\eta)=2|\alpha|^2 +|\mathrm{Trace}(\beta)|^2 .
$
So $\alpha=\mathrm{Trace}(\beta)=0$.
\end{prop}

\begin{proof}
The first term in $(T\eta,\eta)$ is
\begin{equation*}
\begin{split}
(T_1 \eta,\eta)
& :=\sum_{k,\,\ell=1}^N (\rho(X_k)^2 \eta(X_\ell),\eta(X_\ell))_F
=\sum_{k,\,\ell=1}^N (\rho(X_k) \eta(X_\ell),\rho(X_k)\eta(X_\ell))_F \\
& =\sum_{k,\,\ell=1}^N |\beta(X_k,X_\ell)|^2_F
=|\beta|^2 \;\;.
\end{split}
\end{equation*}

The second term in $(T\eta,\eta)$ is
$$
(T_2 \eta,\eta)
:=\sum_{k,\,\ell=1}^N (\rho([X_\ell ,X_k])\eta(X_k),\eta(X_\ell))_F
=(T_3 \eta,\eta)-(T_4 \eta,\eta) \;\;,
$$
where
\begin{equation*}
\begin{split}
(T_3 \eta,\eta)
&:= \sum_{k,\,\ell=1}^N (\rho(X_\ell)\circ \rho(X_k)\eta(X_k),\eta(X_\ell))_F
=\sum_{k,\,\ell=1}^N (\rho(X_k)\eta(X_k),\rho(X_\ell)\eta(X_\ell))_F\\
&= \sum_{k,\,\ell=1}^N (\beta(X_k ,X_k),\beta(X_\ell ,X_\ell))_F
=|\sum_{k=1}^N \beta(X_k ,X_k)|^2_F \\
&=|\mathrm{Trace}(\beta)|^2 \;\;,
\end{split}
\end{equation*}
and
\begin{equation*}
\begin{split}
(T_4 \eta,\eta)
&:= \sum_{k,\,\ell=1}^N (\rho(X_k)\circ \rho(X_\ell)\eta(X_k),\eta(X_\ell))_F
=\sum_{k,\,\ell=1}^N (\rho(X_\ell)\eta(X_k),\rho(X_k)\eta(X_\ell))_F\\
&= (\beta,\beta\circ\phi)\;\;.
\end{split}
\end{equation*}
Here, $\phi \in \textnormal{End}(\LieP \otimes \LieP)$ is defined by
$\phi(X,Y)=(Y,X)$. Note that $\phi$ merely permutes vectors in
the basis of $\mathfrak{p}\otimes\mathfrak{p}$. Therefore
$$
(\sigma,\alpha)=(\sigma\circ\phi,\alpha\circ\phi)=(\sigma,-\alpha)=-(\sigma,\alpha),
$$
thus $(\sigma,\alpha)=0$. Hence
$
|\beta|^2 =|\sigma|^2 +|\alpha|^2
$
and
$$
(\beta,\beta\circ\phi)=
(\sigma+\alpha,\sigma-\alpha)
=|\sigma|^2-|\alpha|^2
=|\beta|^2-2|\alpha|^2 \;\;.
$$
Summing up,
\begin{eqnarray*}
(T\eta,\eta)=|\beta|^2 +|\mathrm{Trace}(\beta)|^2
-(|\beta|^2-2|\alpha|^2)=2|\alpha|^2 +|\mathrm{Trace}(\beta)|^2 .
\end{eqnarray*}
The last assertion follows from Theorem \ref{matu}.
\end{proof}

\section{Consequences of Matsushima-Murakami's method} \label{conseq}

\subsection{Restriction on $S^2\VC^*$-harmonic one-forms}

From now on, $\LL=\SU(n,1)$, $K= U(n)$ and $F = S^2 \VC^*$ is the space of complex
quadratic forms on $\bc^{n+1}$, with the usual action of
$\GL(n+1,\bc)$, $(X,Q)\mapsto X^{\top}QX$, restricted to $\SU(n,1)$. The
admissible inner product on $F$ is the usual $\U(n+1)$-invariant
Hermitian form.

Let
$\mathfrak{su}(n,1)=\mathfrak{u}(n)\oplus \mathfrak{p}$ be the Cartan
decomposition of $\mathfrak{su}(n,1)$.
Here, $\mathfrak{u}(n)=s(\mathfrak{u}(1)\oplus \mathfrak{u}(n))$
consists of traceless block-diagonal skew-Hermitian complex $(n+1)\times
(n+1)$ matrices, and $\mathfrak{p}$ consists of complex matrices of the form
$\left ( \begin{smallmatrix}
0 & x   \\
x^{*}&0
\end{smallmatrix} \right )$, $x\in\bc^n$.

\begin{defi}
We denote by $\chi : U(n) \lo \CC^*$ the standard character $\det$.
\end{defi}

The $SU(n,1)$-module $\VC$
decomposes as a $U(n)$-module~:
$$\VC = \LieP \otimes \chi^{-1} \oplus \chi^{-1} \;\;, $$
(notice that $\LieP \otimes \chi^{-1}$ is nothing else than the
standard $U(n)$-module $\CC^n$).
Thus $S^2 \VC^*$ decomposes as $U(n)$-modules as
$$S^2 \VC^* = (S^2 \LieP^* \oplus \LieP^* \oplus \CC) \otimes \chi^2 \;\;$$
(notice that the $U(n)$-module $S^2 \LieP^*$ is nothing else than
$S^2 \VC^* \cap \LieSp(n)$)
and $\HHom(\LieP, S^2\VC^*)$ as~:
$$ \HHom(\LieP, S^2\VC^*) = (\HHom(\LieP, S^2 \LieP^*) \oplus
\textnormal{End}\, \LieP^* \oplus \LieP) \otimes \chi^2 \;\;.$$
As $U(n)$-modules, $\LieP$ and $S^2 \LieP^*$ are
$\bc$-linear. Thus the $U(n)$-module $\HHom(\LieP, S^2 \LieP^*)$
contains as a direct factor
$\HHom_{\bc}(\LieP, S^2\LieP^{*})$, which contains itself as a
direct factor $S^3 \LieP^*$.

\begin{prop}
\label{pi}
Let $\Gamma$ be a cocompact lattice in $G=SU(n,1)$. Let $\alpha$ be a
$\Gamma$-equivariant harmonic $\ww$-valued 1-form on
$\HypC^n$. Then, for all $Y\in T_{x_0}\HypC^n=\mathfrak{p}$,
$\alpha_{x_0}(Y)\in \wwprime \otimes \chi^2$. Furthermore, $\alpha_{x_0}\in
\HHom_{\br}(\mathfrak{p},S^2\LieP^*)\otimes \chi^2$ is $\bc$-linear and belongs to
the summand $S^3 \LieP^* \otimes \chi^2$ of $\HHom_{\bc}(\LieP,
S^2\LieP^{*}) \otimes \chi^2$.
\end{prop}

\begin{rem}
\label{pibar}
Since $S^2 \bar{V}_{\bc}^*$ is the conjugate vector space of $S^2 V_{\bc}^*$, Proposition \ref{pi} implies that $\Gamma$-equivariant harmonic $S^2 \bar{V}_{\bc}^*$-valued 1-forms on $H^n_\bc$ are in fact $S^2 \bar{\mathfrak{p}}^{*}\otimes \chi^2$-valued $(0,1)$-forms.
\end{rem}

\subsection{Proof of Proposition \ref{pi}}

 A straightforward calculation yields

\begin{lem}
\label{adjoint}
Let $X=\left (\begin{smallmatrix}
0 & x   \\
x^{*} & 0
\end{smallmatrix} \right )$, $x\in\bc^n$, be a vector of
$\mathfrak{p}$. Let $Q=\left (\begin{smallmatrix}
A & B   \\
B^{\bot} & d
\end{smallmatrix} \right ) \in S^2 \VC^*$. Then
\begin{eqnarray*}
\rho(X)(Q)=X^{\top}Q+QX=\left (  \begin{smallmatrix}
(Bx^{*})^{\top}+Bx^{*} & Ax+d\bar{x}\\
(Ax+d\bar{x})^{\top} &  2x^{\top}B
\end{smallmatrix} \right ).
\end{eqnarray*}
\end{lem}

Let $\eta\in \HHom(\mathfrak{p},S^2\VC^*)$ be represented by a matrix
$Q=Y\mapsto \left ( \begin{smallmatrix}
A(Y) & B(Y)   \\
B(Y)^{\bot} & d(Y)
\end{smallmatrix} \right )$ of $\br$-linear forms on
$\mathfrak{p}$. Then the bilinear form $\beta=\left ( \begin{smallmatrix}
\beta_{1}&\beta_{2}\\
\beta_{2}^{\top} & \beta_{3}
\end{smallmatrix} \right )$
becomes a triple of matrix valued bilinear forms on $\bc^n$,
\begin{eqnarray*}
\beta_{1}(x,y)&=&(B(y)x^{*})^{\top}+B(y)x^{*},\\
\beta_{2}(x,y)&=&A(y)x+d(y)\bar{x},\\
\beta_{3}(x,y)&=&2x^{\top}B(y).
\end{eqnarray*}

According to Proposition \ref{calculT}, $(T\eta,\eta)=0$ if and only if the following 6 equations hold.
\begin{eqnarray*}
\beta_{1},~\beta_{2}~\mathrm{and}~\beta_{3}~\mathrm{are~symmetric},\;\;
\mathrm{Trace}(\beta_{1})=0, \;\;\mathrm{Trace}(\beta_{2})=0, \;\;
\mathrm{Trace}(\beta_{3})=0.
\end{eqnarray*}

\begin{lem}
\label{beta3}
If $\beta_{3}$ and $\beta_1$ are symmetric, then $B=0$.
\end{lem}

\begin{proof}
Let $B_{\bc}$ and $B_{\bar{\bc}}$ denote the $\bc$-linear (resp. anti
$\bc$-linear) components of the $\br$-linear map $B$. Matrixwise, each
of $B_{\bc}$ and $B_{\bar{\bc}}$ is given by a $n\times n$ complex
matrix $\mathcal{B}_{\bc}$ (resp. $\mathcal{B}_{\bar{\bc}}$), and
$B(y)=\mathcal{B}_{\bc}y+\mathcal{B}_{\bar{\bc}}\bar{y}$. Thus
\begin{eqnarray*}
\beta_3 (x,y)=x^{\top}\mathcal{B}_{\bc}y+x^{\top}\mathcal{B}_{\bar{\bc}}\bar{y}
\end{eqnarray*}
is the sum of a $\bc$-bilinear and a sesquilinear form. If $\beta_3$
is symmetric, the sesquilinear part vanishes
(i.e. $\mathcal{B}_{\bar{\bc}}=0$), and $\mathcal{B}_{\bc}$ is
symmetric.

Next,
\begin{eqnarray*}
\beta_1 (x,y)=(\mathcal{B}_{\bc}yx^{*})^{\top}+\mathcal{B}_{\bc}yx^{*}
\end{eqnarray*}
is sesquilinear. If $\beta_1$ is symmetric, it is identically
zero. Since rank one matrices of the form $yx^{*}$ span all $n\times
n$ complex matrices,
$(\mathcal{B}_{\bc}M)^{\top}+\mathcal{B}_{\bc}M=0$ for all $n\times n$
complex matrices $M$. Take $M=\mathcal{B}_{\bc}^{*}$ and take the
trace to conclude that $\mathcal{B}_{\bc}=0$.
\end{proof}

\begin{lem}
\label{beta2}
If $\beta_{2}$ is symmetric and $\mathrm{Trace}(\beta_{2})=0$, then
$d=0$ and $A(y)$ depends $\bc$-linearly on $y$. Furthermore,
identifying $\bc^n$-valued bilinear maps with trilinear forms,
$(x,y)\mapsto A(y)x$ is fully symmetric.
\end{lem}

\begin{proof}
Let $A_{\bc}$ and $A_{\bar{\bc}}$ denote the $\bc$-linear (resp. anti
$\bc$-linear) components of the $\br$-linear map $A:\bc^n \ra
S^{2}(\bc^n)$. Similarly, let $d_{\bc}$ and $d_{\bar{\bc}}$ denote the
$\bc$-linear (resp. anti $\bc$-linear) components of the $\br$-linear
form $d$. If $\beta_2 : (x,y)\mapsto A(y)x+d(y)\bar{x}$ is symmetric,
then
\begin{eqnarray*}
\begin{cases}
\forall x,\,y\in\bc^n ,~A_{\bar{\bc}}(y)x &=d_{\bc}(x)\bar{y},\\
(x,y)\mapsto d_{\bar{\bc}}(y)\bar{x}  & \text{is symmetric},\\
(x,y)\mapsto A_{\bc}(y)x  & \text{is symmetric}.
\end{cases}
\end{eqnarray*}

The trace of the restriction of $\beta$ to a complex line $\bc e$,
$|e|=1$, depends only on its sesquilinear part
\begin{eqnarray*}
\beta_2^{sq}(x,y)=A_{\bar{\bc}}(y)x+d_{\bc}(y)\bar{x}=d_{\bc}(x)\bar{y}+d_{\bc}(y)\bar{x}.
\end{eqnarray*}
and is equal to $2\beta_2^{sq}(e,e)=4d_{\bc}(e)\bar{e}$. Let $e_1
,\ldots,e_n$ be a Hermitian basis of $\bc^n$. Then
\begin{eqnarray*}
\mathrm{Trace}(\beta_2)=\mathrm{Trace}(\beta_2^{sq})=4\sum_{k=1}^{n}d_{\bc}(e_k)\bar{e}_k .
\end{eqnarray*}
Since $\mathrm{Trace}(\beta_{2})=0$, we get $d_{\bc}=0$. This implies
that $A_{\bar{\bc}}(y)x=0$ for all $x$ and $y$,
i.e. $A_{\bar{\bc}}=0$.

Next, pick a nonzero vector $y\in\mathrm{ker}(d_{\bar{\bc}})$. Since
$(x,y)\mapsto d_{\bar{\bc}}(y)\bar{x}$ is symmetric, for all
$x\in\bc^n$, $d_{\bar{\bc}}(x)\bar{y}=0$, thus $d_{\bar{\bc}}=0$.

Finally, view the components of $A_{\bc}(y)x$ in some Hermitian
basis $e_1,\cdots,e_n$ of  $\bc^n$ as bilinear forms on $\bc^n$,
with respective matrices
$\mathcal{A}^{1}=A_\bc(e_1),\ldots,\mathcal{A}^{n}=A_\bc(e_n)$.
Since the values $A_\bc(y)$ are symmetric matrices, these matrices
are symmetric, $\mathcal{A}^{\ell}_{jk}=\mathcal{A}^{\ell}_{kj}$.
But for every $y=(y_{1},\ldots,y_{n})\in\bc^n$,
$$(A_\bc(y)x)_j=y_{\ell}\mathcal{A}^{\ell}_{jk}x_k=(A_\bc(x)y)_j=x_k\mathcal{A}^k_{j\ell}y_\ell.$$
This implies that
 $\mathcal{A}^{\ell}_{jk}=\mathcal{A}^{k}_{j\ell}$. Hence $\mathcal{A}^{\ell}_{jk}$ is fully symmetric.
\end{proof}

\section{Second order obstruction} \label{second}

\subsection{Cup-product, case $\G=\Sp(n,1)$}

\begin{defi}
\label{deflambda}
Let $\lambda: \LieSp(n,1) \lo \RR$ be the $SU(n,1)$-invariant linear
form defined by the Killing inner product with the
$SU(n,1)$-invariant vector $iI_{n+1}$, which generates the centralizer
of $SU(n,1)$ in $Sp(n,1)$.
\end{defi}
The restriction of the Killing form of $Sp(n,1)$ to $Sp(n)$ is
proportional to the Killing form of $Sp(n)$, which is proportional to
$\Re e(\mathrm{Trace}_{\bh}(A^{*}A))$. Therefore, for $A\in
\mathfrak{sp}(n)\subset \mathfrak{sp}(n,1)$,
\begin{eqnarray*}
\lambda(A)=A\cdot iI_{n+1}=-\Re e(i\,\mathrm{Trace}_{\bh}(A)),
\end{eqnarray*}

\begin{lem}
\label{square}
Let $\alpha$ be an $\mathfrak{sp}(n)$-valued $(1,0)$-form
on $T_{x_{0}}\HypC^n=\mathfrak{p}$. Assume that $\alpha$ belongs to
$Hom_{\bc}(\mathfrak{p},\wwprime)$. Thanks to the Lie bracket
of $\mathfrak{sp}(n)$, $[\alpha, \alpha]$ becomes an
$\mathfrak{sp}(n)$-valued 2-form on $\mathfrak{p}$. Let $\omega$
denote the K\"ahler form on $\mathfrak{p}$. There is a nonzero
constant $c$ such that
\begin{eqnarray*}
\lambda\circ[\alpha, \alpha]\wedge\omega^{n-1}=c|\alpha|^2
\omega^{n}.
\end{eqnarray*}
\end{lem}

\begin{proof}
Recall that the embedding of $\wwprime$ to $\mathfrak{sp}(n)$ is
  defined by
$A\mapsto jQA$ where $Q=(I_n,-1)$ a diagonal matrix.
Write $\alpha=jQ\delta$
                     where $\delta$ is a
symmetric complex matrix of $(1,0)$-forms. Then, for all $Y$,
$Y'\in\mathfrak{p}$,
\begin{eqnarray*}
\alpha\wedge\alpha(Y,Y')=
\alpha(Y)\otimes\alpha(Y')-\alpha(Y')\otimes\alpha(Y)\in\mathfrak{sp}(n)\otimes\mathfrak{sp}(n)\;\;,
\end{eqnarray*}
\begin{eqnarray*}
[\alpha, \alpha](Y,Y')=[\alpha(Y),\alpha(Y')]-[\alpha(Y'),\alpha(Y)]=2[\alpha(Y),\alpha(Y')]\in\mathfrak{sp}(n).
\end{eqnarray*}
Let $A$, $B$ be two symmetric complex matrices. The Lie bracket of their images in $\mathfrak{sp}(n)$ is
\begin{eqnarray*}   [jQA,jQB]=jQAjQB-jQBjQA
=-\bar{A}B+\bar{B}A,
\end{eqnarray*}
(note it belongs to $\mathfrak{u}(n)$), thus
\begin{eqnarray*}
[\alpha, \alpha](Y,Y')=-2(\overline{\delta(Y)}\delta(Y')-\overline{\delta(Y')}\delta(Y)),
\end{eqnarray*}
showing that $[\alpha, \alpha]$ is a matrix of $(1,1)$-forms.
Up to a nonzero constant,
\begin{eqnarray*}
\lambda\circ[\alpha,\alpha](Y,Y')=\Im
m(\mathrm{Trace}_{\bc}(\delta(Y)^{*}\delta(Y'))).
\end{eqnarray*}
Note that
$\lambda\circ[\alpha,\alpha](Y,iY)=\mathrm{Trace}(\delta(Y)^{*}\delta(Y))=|\delta(Y)|^2
>0$.

If $\phi$ is a $(1,1)$-form on $\bc^n$, then
\begin{eqnarray*}
\frac{\phi\wedge\omega^{n-1}}{\omega^{n}}=\frac{2}{n}\sum_{k=1}^{n}\phi(E_k ,iE_k),
\end{eqnarray*}
where $E_1 ,\ldots,E_n$ is a unitary basis of $\bc^n$ (i.e. $(E_1 , iE_1 ,\ldots,E_n ,iE_n)$ is an orthonormal
basis of the underlying real Euclidean vectorspace). Therefore
\begin{eqnarray*}
\frac{\lambda\circ[\alpha,\alpha]\wedge\omega^{n-1}}{\omega^{n}} =
\frac{2}{n}\sum_{k=1}^{n}|\delta(E_k)|^2
\end{eqnarray*}
is a nonzero multiple of $|\alpha|^2$.
\end{proof}

The following proposition finishes the proof of
theorem~\ref{main_theorem}, in case $\G=\Sp(n,1)$:
\begin{prop}
\label{cup}
Let $\alpha$ be a nonzero harmonic
$\mathfrak{sp}(n,1)$-valued 1-form on $\Gamma \backslash \HypC^n$. Assume that the
component of $\alpha$ on the centralizer of $\mathfrak{su}(n,1)$ in
$\mathfrak{sp}(n,1)$ vanishes. Then $[\alpha, \alpha]\not=0$ in
$H^2 (\Gamma,\mathfrak{sp}(n,1))$. In particular, $\alpha$ does not
integrate into a nontrivial deformation of the conjugacy class of the
embedding $\Gamma\hookrightarrow Sp(n,1)$.
\end{prop}

\begin{proof}
By contradiction. According to Weil's vanishing theorem, the
$\mathfrak{su}(n,1)$-component of $\alpha$ vanishes, thus $\alpha$ is
$\ww$-valued. According to Proposition \ref{pi}, $\alpha$ can
be viewed as a smooth section of the homogeneous bundle over $\Gamma
\backslash \HypC^n$ whose
fiber is the subspace $S^3 \LieP^*$ of
$\HHom_{\bar{\bc}}(\mathfrak{p},\wwprime)$. In particular, $\alpha$
can be viewed pointwise as a $\mathfrak{sp}(n)$-valued
$(1,0)$-form. Assume that $[\alpha, \alpha]=0$, i.e., that there
exists a $\mathfrak{sp}(n,1)$-valued 1-form
$\eta$ on $M$
such that $d\eta=[\alpha, \alpha]$. Then, with Lemma \ref{square},
\begin{equation*}
\begin{split}
2c\parallel \alpha \parallel_{L^{2}(M)}^{2}
&= \int_{M}c |\alpha|^{2}\omega^{n}
=\int_{M}\lambda\circ[\alpha, \alpha]\wedge\omega^{n-1}\\
&=\int_{M}\lambda\circ(d\eta)\wedge\omega^{n-1}
=\int_{M}d(\lambda\circ\eta\wedge\omega^{n-1})
=0\;\;,
\end{split}
\end{equation*}
thus $\alpha=0$, contradiction.
\end{proof}

\subsection{Cup-product, case $\G=\SO(4n,4)$}

Choose $J=\left (  \begin{smallmatrix}
0 & -I_{n+1}\\
I_{n+1} & 0
\end{smallmatrix} \right )$ as a complex structure on $\br^{2n+2}$.

\begin{defi}
\label{deflambdapp}
On $\mathfrak{so}(4n,4)$, there is a $\SU(n,1)$-invariant linear forms
$\lambda'$ and $\lambda''$, given by Killing inner product with the $\SU(n,1)$-invariant vectors $J'=\left (  \begin{smallmatrix}J&0\\ 0&-J\end{smallmatrix} \right )$ and $J''=\left (  \begin{smallmatrix}J&0\\ 0&J\end{smallmatrix} \right )$.
\end{defi}

\begin{prop}
\label{lambdappp}
$\lambda'$ vanishes on
$[S^2 \bar{V}_{\bc}^*,S^2 \bar{V}_{\bc}^*]$, $\lambda''$ vanishes on
$[S^2 V_{\bc}^*,S^2 V_{\bc}^*]$, and they both vanish on $[S^2
V_{\bc}^*, S^2 \bar{V}_{\bc}^*]$.
\end{prop}

\begin{proof}
Let $Z=\left(\begin{smallmatrix}0&BQ'\\ -B^*Q'&0\end{smallmatrix}\right)$,
$Z'=\left(\begin{smallmatrix}0&B'Q'\\ -B'^* Q'&0\end{smallmatrix}\right)\in S^2
V_{\bc}^*\oplus S^2 \bar{V}_{\bc}^*$. View $Z$ and $Z'$ as vectors
in $\mathfrak{so}(4n,4)$. Then using the fact that $Q'$ commutes
with $B,B',B^*,B'^*$,
\begin{eqnarray*}
\lambda'([Z,Z'])&=&\mathrm{Trace}(J[-BB'^*+B'B^*+B^*B'-B'^*B])\\
\lambda''([Z,Z'])&=&\mathrm{Trace}(J[-BB'^*+B'B^*-B^*B'+B'^*B]).
\end{eqnarray*}
If both $Z$ and $Z'$ are in $S^2 V_{\bc}^*$,
$$\lambda'([Z,Z'])=2\mathrm{Trace}(J[B'B^*-BB'^*])$$
and vanishes if both are in $S^2 \bar{V}_{\bc}^*$, using the fact
that $J$ anti-commutes with $B'$ and $B'^*$.

If both $Z$ and $Z'$ are in $S^2 \bar{V}_{\bc}^*$,
$$\lambda''([Z,Z'])=2\mathrm{Trace}([B'B^*-BB'^*]J)$$
and vanishes if both are in $S^2 V_{\bc}^*$.

If $Z\in S^2 V_{\bc}^*$ and $Z'\in S^2 \bar{V}_{\bc}^*$, both $\lambda'$, $\lambda''$ vanish. Indeed, since $BJ=JB$ and $B'J=-JB'$, for example,
$\mathrm{Trace}(JB'B^*)=\mathrm{Trace}(B'B^*J)=\mathrm{Trace}(-B'JB^*)=\mathrm{Trace}(-B'B^*J)=0$.
\end{proof}

If $Z$ and $Z'\in S^2 \mathfrak{p}^* \subset S^2 V_{\bc}^*$, write
$B=\left(\begin{smallmatrix}C&-D\\ D&C\end{smallmatrix}\right)$ and
$B'=\left(\begin{smallmatrix}C'&-D'\\ D'&C'\end{smallmatrix}\right)$ where $C$, $C'$,
$D$, $D'$ are symmetric.
\begin{eqnarray*}
\lambda'([Z,Z'])=8\mathrm{Trace}(DC'-CD').
\end{eqnarray*}
The complex structure on $S^2 V_{\bc}^*$ is $B\mapsto \mathcal{J}(B)=JB$, i.e. $(C,D)\mapsto (-D,C)$. Thus
\begin{eqnarray*}
\lambda'([Z,\mathcal{J}(Z)])&=&8\mathrm{Trace}(-D^{2}-C^{2})\\
&=&-8\mathrm{Trace}(D^{\top}D+C^{\top}C)=-4|B|^2 =-2|Z|^2 .
\end{eqnarray*}

If $Z$ and $Z'\in S^2 \bar{\mathfrak{p}}^* \subset S^2
\bar{V}_{\bc}^*$, write $B=B^*=\left(\begin{smallmatrix}C&D\\
D&-C\end{smallmatrix}\right)$ and $B'=B'^*=\left(\begin{smallmatrix}C'&D'\\
D'&-C'\end{smallmatrix}\right)$ where $C$, $C'$, $D$ and $D'$ are symmetric.
Then
\begin{eqnarray*}
\lambda'([Z,Z'])=8\mathrm{Trace}(DC'-CD').
\end{eqnarray*}
The complex structure on $S^2 \bar{V}_{\bc}^*$ is $B\mapsto \mathcal{J}(B)=-JB$, i.e. $(C,D)\mapsto (D,-C)$. Thus
\begin{eqnarray*}
\lambda'([Z,\mathcal{J}(Z)])&=&8\mathrm{Trace}(D^{2}+C^{2})\\
&=&8\mathrm{Trace}(D^{\top}D+C^{\top}C)=4|B|^2 =2|Z|^2 .
\end{eqnarray*}

Let $\eta$ be an equivariant harmonic $\mathfrak{so}(4n,4)$-valued
1-form. According to \cite{rag} and Proposition \ref{pi},
$\eta=\tau+\alpha+\alpha'$ where $\tau$ is $\mathfrak{z}$-valued,
$\alpha$ is a $S^2 \mathfrak{p}^*$-valued $(1,0)$-form and $\alpha'$
a $S^2 \bar{\mathfrak{p}}^*$-valued $(0,1)$-form.

We have seen that
$$\lambda'\circ([\alpha\wedge\alpha'])=0,$$
$$\lambda''\circ([\alpha\wedge\alpha'])=0.$$
 Therefore
$$
(\lambda'+\lambda'')\circ([(\alpha+\alpha')\wedge(\alpha+\alpha')])=\lambda'\circ([\alpha\wedge\alpha])+\lambda''\circ([\alpha'\wedge\alpha']).
$$
If $Y\in\mathfrak{p}$, since $\alpha$ has type $(1,0)$, $\alpha(iY)=\mathcal{J}(\alpha(Y))$,
\begin{eqnarray*}
\lambda'\circ([\alpha\wedge\alpha])(Y,iY)=2|\alpha(Y)|^2 .
\end{eqnarray*}

Since $\alpha'$ has type $(0,1)$, $\alpha'(iY)=-\mathcal{J}(\alpha'(Y))$,
\begin{eqnarray*}
\lambda''\circ([\alpha'\wedge\alpha'])(Y,iY)=2|\alpha'(Y)|^2 .
\end{eqnarray*}

It follows that
\begin{eqnarray*}
\frac{(\lambda'+\lambda'')\circ[(\alpha+\alpha')\wedge(\alpha+\alpha')]\wedge\omega^{n-1}}{\omega^{n}}=\frac{2}{n}(|\alpha|^2
+|\alpha'|^2).
\end{eqnarray*}
Again, if the cohomology class of
$[(\alpha+\alpha')\wedge(\alpha+\alpha')]$ vanishes, then the $L^2$
norm of $\alpha$ and $\alpha'$ vanishes. This shows that the
quadratic map induced by bracket-cup product on
$H^{1}(\Gamma,\mathfrak{so}(4n,4))/H^{1}(\Gamma,\mathfrak{z})$ is
anisotropic.

\section{A more geometric proof}  \label{nonab}

In this section, we sketch a second proof of the proposition~\ref{pi}
for $\G = \Sp(n,1)$, using a theorem of Carlson-Toledo \cite{cto} and
some non-Abelian Hodge theory.

\subsection{Reminder on quaternionic K\"ahler manifolds} \label{quaternionic_geometry}
For the convenience of the reader we recall some
general facts on quaternionic K\"ahler manifolds. We refer to
\cite{sal} for a panorama.

\begin{defi}
A Riemannian manifold $M$ of dimension $4n$ is quaternionic K\"ahler if its holonomy
group is contained in the subgroup $Sp(n) Sp(1):= Sp(n) \times_{\ZZ/
  2\ZZ} Sp(1)$ of $SO(4n)$, where $\ZZ/2\ZZ$ is generated by $(-I, -1)$.
\end{defi}

It is well known that such a Riemannian manifold $M$ is always Einstein.
Even if $M$ is not necessarily K\"ahler, its
geometry can be essentially understood from the point of view of
complex geometry.
\begin{defi} We denote by $\mathcal{P}_M$ the canonical $Sp(n,1)
Sp(1)$-reduction of the principal bundle of
orthogonal frames of $M$, and by $\E_M$ the canonical $3$-dimensional
parallel sub-bundle $\mathcal{P}_M \times_{Sp(n) Sp(1)} \RR^3$ of
$\End(TM)$.
\end{defi}

\begin{defi}
Let $p: Z \lo M$ be the $S^2$-fiber bundle on $M$ associated to the action of
$Sp(1)/\ZZ_2 \simeq SO(3)$ on $S^2$~: $Z = \mathcal{P}_M \times_{Sp(n)
  Sp(1)} S^2$. The space $Z$ is called the
twistor space of $M$.
\end{defi}

In other words, $Z$ is the unit sphere of $\E_M$.

\begin{theor} \cite[theor. 4.1]{sal0}
Let $M$ be a quaternionic K\"ahler manifold. Then its twistor space $Z$
admits a canonical complex structure, for which the fibers of $p:Z \lo
M$ are complex rational curves.
\end{theor}

As $\E_M$ is parallel, it inherits from the Levi-Civita connection on $TM$
a linear connection compatible with the metric. It follows that the
corresponding horizontal distribution induces a horizontal distribution
$T^\RR_h(Z) \subset T^\RR(Z)$. In the case where the scalar curvature of $M$ is
non-zero, one can show that $T^\RR_h(Z)$ naturally defines an
horizontal holomorphic distribution $T_h Z \subset TZ$ making $Z$ a
holomorphic contact manifold \cite[prop. 5.2]{sal}.

\sspace
An other ingredient of some importance for us is the following~:
\begin{lem}
Any quaternionic K\"ahler manifold $M$ admits a non-zero closed
$4$-form $\Omega_M$, canonical up to homothety.
\end{lem}
\begin{proof}
Just notice that the $Sp(n)Sp(1)$-module $\bigwedge^4 (\RR^{4n})^*$ admits
a unique trivial submodule of rank $1$.
\end{proof}

\begin{lem} \label{pontryagin}
The form $\Omega_M$ (conveniently normalized) is the Chern-Weil form of the first Pontryagin
class $p_1(\E_M) \in H^4(M , \ZZ)$.
\end{lem}

\begin{proof}
This is proved in \cite[p.148-151]{sal0}.
\end{proof}

\subsection{Quaternionic K\"ahler symmetric spaces}
The description of quaternionic K\"ahler symmetric spaces and their
twistor spaces is due to Wolf \cite{wol}, following Boothby
\cite{boo}. There exists $3$ families of quaternionic K\"ahler non-compact
irreducible symmetric spaces of dimension $4n$, $n \geq 2$~: $\HypH^n = Sp(n,1)/Sp(n) \cdot Sp(1)$,
$X^n= SU(n,2)/S(U(n) \times U(2))$ and $Y^n =SO(n,4)/S(O(n) \times O(4))$.
The only K\"ahler one is $X^n$. In each case the isotropy group is of
the form $K\cdot Sp(1)$ and the twistor space is obtained by
replacing the $Sp(1)$-factor by  $U(1)$. Notice that the twistor
map for $X^n$ is {\em not holomorphic}.

By functoriality of the twistor construction, we associate to the sequence of
totally geodesic quaternionic K\"ahler embeddings $\HypH^n \hookrightarrow X^{2n} \lo Y^{4n}$
the
commutative diagram~:
\begin{equation} \label{le_diagramme}
\xymatrix{
D^{\HH}_{2n+1}= \frac{Sp(n,1)}{Sp(n) \times U(1)} \ar[d]_p \ar@{^(->}[r] &D^\CC_{4n+1}=
\frac{SU(2n,2)}{S(U(1) \times U(2n) \times U(1))} \ar[d]_p \ar@{^(->}[r]
&D^\RR_{8n+1}= \frac{SO(4n,4)}{S(O(4n) \times U(2))}
\ar[d]_p \\
\HypH^n \ar@{^(->}[r] &X^{2n} \ar@{^(->}[r] &Y^{4n}
}\;\;,
\end{equation}
where the vertical maps are twistor fibrations and the horizontal maps on
the top line are holomorphic closed {\em horizontal} (i.e. preserving the contact
structure) immersions.

\subsection{An invariant for quaternionic representations}

Let $X$ a smooth manifold and $\rho: \Ga= \pi_1(X) \lo Sp(n,1)$ a
representation. Choose any $\rho$-equivariant smooth map $\phi: \tilde{X} \lo
\HypH^n$. The pull-back $\phi^* \E_{\HypH^n}$ is a
$\Ga$-equivariant rank $3$ real bundle on $\tilde{X}$. Thus it descends
to a bundle on $X$, still denoted $\phi^* \E_{\HypH^n}$. By
lemma~\ref{pontryagin} and the functoriality of characteristic classes, the
$4$-form $\phi^* \Omega_{\HypH^n}$ represents the Pontryagin class $p_1(\phi^*
\E_{\HypH^n}) \in H^4(X, \ZZ)$. As $\HypH^n$ is a
contractible space, any two $\rho$-equivariant
maps $\phi, \phi': \tilde{X} \lo \HypH^n$ are $\rho$-equivariantly
homotopic. Finally the class $[\phi^* \Omega_{\HypH^n}]\in H^4(X,
\ZZ)$ depends only on $\rho$.

\begin{defi}
Let $X$ a smooth manifold and $\rho: \Ga= \pi_1(X) \lo Sp(n,1)$ a
representation. We denote by $c_{\rho}\in H^4(X, \ZZ)$ the class $[\phi^* \Omega_{\HypH^n}] \in H^4(X, \ZZ)$.
\end{defi}

\begin{rem}
The invariant $c_\rho$ is a quaternionic version of the (Hermitian)
Toledo invariant.
\end{rem}

\begin{lem} \label{constant}
Let $M$ be a smooth manifold and $\Ga = \pi_1(M)$.
The function $$c: \M(\Ga, \Sp(n,1))(\RR) \lo H^4(X, \ZZ)$$ which
to $[\rho]$ associates $c_\rho$ is constant on connected components of
$\M(\Ga, \Sp(n,1))(\RR)$.
\end{lem}

\begin{proof}
This follows immediately from the integrality of $c_{\rho}$.
\end{proof}

\subsection{Link with Hodge theory}
Twistor spaces of quaternionic K\"ahler symmetric spaces are the
simplest examples of Griffiths's period domains for variations of Hodge
structures. One easily proves
the following lemma (a proof for the global embedding $\HypH^n
\hookrightarrow Y^{4n}$ is provided in \cite[p.192-193]{cto}, the
proof for the other maps is similar)~:

\begin{lem}
Let $K$ be $\RR$, $\CC$ or $\HH$. Let $r_K$ be $4$, $2$ or $1$
respectively.
\begin{itemize}
\item
Each twistor space $D^K_{2r_K \cdot (n +1)}$ is the Griffiths's period
domain for polarized weight $2$ pure Hodge structures with Hodge numbers $(2,
4n, 2)$ on $\RR^{4n+4}$, stable under $K$-multiplication (when we
identify $\RR^{4n+4}$ with $K^{r_{K} \cdot (n+1)}$).
\item The inclusions in the sequence  $D^{\HH}_{2n+1} \hookrightarrow D^\CC_{4n+1}
\hookrightarrow D^\RR_{8n+1}$ correspond to the functors partially
forgetting the $K$-stability condition, for the inclusions $\RR
\subset \CC \subset \HH$.
\end{itemize}
\end{lem}

\subsection{Any deformation is a complex variation of Hodge structures}
Let $n>1$, $\SU(n,1) = \SU(\VC, \hc)$ and let $j: \SU(n,1)
\hookrightarrow \U(n,1)
\hookrightarrow \Sp(n,1)= \SU(\VH, \hh)$ be the natural
embedding. Let $f: \HypC^n \hookrightarrow \HypH^n$ be the corresponding totally
geodesic $U(n,1)$-equivariant embedding. Notice that it canonically
lifts to a holomorphic $U(n,1)$-equivariant embedding $\tilde{f}:
\HypC^n \hookrightarrow  D^\HH_{2n+1}$ making the
$U(n,1)$-equivariant diagram
$$
\xymatrix{
\HypC^n \ar@{^(->}[r]^{\tilde{f}}\ar@{^(->}[rd]_{f} & D^\HH_{2n+1} \ar[d]^{p} \\
&\HypH^n}
\;\;,
$$
commutative.

Let $i: \Gamma \hookrightarrow SU(n,1) = \SU(n,1)(\RR)$ be a cocompact
torsion-free lattice and $\rho = j \circ i : \Ga \lo Sp(n,1)$ the
corresponding representation. Let $M$ be the compact K\"ahler manifold
$\Ga \backslash \HypC^n$.

\begin{lem} \label{nonvanish}
$c_\rho \not = 0 \in H^4(M, \ZZ)/\text{torsion}$.
\end{lem}

\begin{proof}
By lemma~\ref{pontryagin} the (descent to $M$ of
the) curvature form ${f}^* \Omega_{\HypH^n}$
is nothing else than the Chern-Weil form of (the descent to $M$ of)
the bundle
$p_1({f}^*\E_{\HypH^n})$. As the standard representation $\RR^3$ of $SO(3)$
decomposes as $\RR \oplus \CC$ as an $U(1) \subset SO(3)$-module, where
$\RR$ is the trivial representation and $\CC$ is the standard
$U(1)$-module, we obtain that ${f}^*\E_{\HypH^n}$ is the
direct sum of the trivial bundle and (the descent to $M$ of) the holomorphic line bundle $\Lcal = U(n,1)
  \times_{U(n) \times U(1)} \CC \lo \HypC^n$. Thus
$$p_1({f}^*\E_{\HypH^n}) = -c_2(\Lcal \otimes_\RR \CC) =
-c_2(\Lcal \oplus \overline{\Lcal}) = c_1^2(\Lcal)\;\;.$$
Finally
${f}^*\Omega_{\HypH^n} =  \omega_{M}^2 $, where $\omega_{M}$ is the standard
    K\"ahler form on the quotient $M$ of $\HypC^n$.
Thus $c_{\rho} = [\omega_M]^2 \not = 0 \in H^4(M, \ZZ)/\text{torsion}$
(as $n>1$).
\end{proof}

Let $\lambda : \Ga \lo Sp(n,1)$ be a reductive representation in the same
connected component of $M(\Ga, \Sp(n,1))(\RR)$ as $\rho$. By
lemma~\ref{constant} and ~\ref{nonvanish}, $c_\lambda = c_\rho \not =0
\in H^4(M, \ZZ)/\text{torsion}$. Let $f_\lambda : \tilde{M} \lo
\HypH^n$ be the $\lambda$-equivariant harmonic map. As $c_{\lambda} =[f_\lambda^* \Omega_{\HypH^n}] \in H^4(M, \RR)$, the
harmonic map $f_\lambda$ is of rank at least $4$ on some open
subset of $\tilde{M}$. Thus we can apply the following result of
Carlson-Toledo~:
\begin{theor} \cite[theor. 6.1.]{cto} \label{carto}
Let $X$ be a compact K\"ahler manifold, $\lambda : \Ga= \pi_1(X) \lo
\Sp(n,1)(\RR)$ a reductive representation and $f_\lambda: \tilde{X} \lo
\HypH^n$ the $\rho$-equivariant harmonic map, where $\tilde{X}$ denotes the
universal covering of $X$. {\em Assume that the rank of the
  differential $df: T\tilde{X} \lo T\HypH^n$ is larger than $2$ at
  some point $x$ of $X$}. Then there exists a horizontal holomorphic
$\lambda$-equivariant period map
$\tilde{f}_\lambda : \tilde{X} \lo D^{\HH}_{2n+1}$ making the following diagram commute~:
$$
\xymatrix{
\tilde{X} \ar@{^(->}[r]^{\tilde{f}_\lambda} \ar@{^(->}[rd]_{f_\lambda} & D^\HH_{2n+1} \ar[d]^{p} \\
&\HypH^n}
\;\;.
$$
\end{theor}
Thus we obtain that any deformation $\lambda$ of $\rho$ is still the monodromy of a variation of
Hodge structure $\tilde{f}_\lambda: \tilde{M} \lo D^{\HH}_{2n+1}$. To
prove proposition~\ref{pi} is thus equivalent to the following~:

\begin{prop} \label{pi'}
The tangent space at $(f, \rho)$ to the space of $\Sp(n,1)$-variations
of Hodge structures identifies (as a real vector space) with~:
$$H^1(M, \LieZ_{\mathfrak{u}(n,1)}\mathfrak{su}(n,1)) \oplus H^0(M, S^3 T^*M \otimes
L_{\chi^2})\;\;,$$
where $\LieZ_{\mathfrak{u}(n,1)}\mathfrak{su}(n,1) \simeq \RR$ denotes
the Lie algebra of the centralizer $S^1$ of $SU(n,1)$ in $U(n,1)$ and
$L_{\chi^2} = \Gamma \backslash SU(n,1) \times_{U(n), \chi^2} \CC$ denotes the automorphic line bundle
on $M$ associated to the character $\chi^2 : U(n) \lo S^1$.
\end{prop}

\subsection{Proof of proposition~\ref{pi'}}

\subsubsection{Explicit notations}
We fix a basis $(e_0, \cdots, e_n)$ of $\VR$ over $\RR$. This is also
a basis of $\VC$ over $\CC$. As a $(2n+2)$-$\CC$-vector space, $\VH$ is isomorphic to $\VC \oplus j
\VC$ and we choose $(f_0= e_o, \cdots, f_n=e_n, f_{n+1}=-e_n \cdot j
, \cdots, f_{2n}= -e_1 \cdot j, f_{2n+1}=
e_0 \cdot j)$.
The right multiplication by $j \in \HH$ on a column vector $\mathbf{v}
= ( v_0,  \cdots, v_{2n})^t $ in the basis $(f_i)_{0 \leq i \leq 2n+1}$ is given by
$ \mathbf{v} \cdot j=
\left (
\begin{smallmatrix}
& -J_1 \\
J_1^t &
\end{smallmatrix} \right ) \cdot \overline{\mathbf{v}}$, where
$J_1$ denotes the $(n+1) \times(n+1)$-matrix
$$
J_1 =  \left (
\begin{smallmatrix}
& & &-1 \\
& &1 & \\
& \adots & & \\\
1 & & &
\end{smallmatrix} \right ) \;\;.
$$
This realizes the group $\GL(n+1, \HH)$ of $\HH$-linear automorphism of
$\VH$ as the
matrix group
$$ \GL(n+1,\HH ) = \{ X \in \GL(2n+2, \CC) \; / \; X \cdot
\left ( \begin{smallmatrix}
& J_1 \\
-J_1^t &
\end{smallmatrix} \right )
=
\left ( \begin{smallmatrix}
& J_1 \\
-J_1^t &
\end{smallmatrix} \right )
\cdot \overline{X} \} \;\;.$$

The real orthogonal form $Q_{\RR}$ of signature $(n,1)$
on $\VR$ is defined by~:
$Q_{\RR}(x_0 \cdot e_0 + \cdots x_n e_n )= -x_0^2 + x_1^2 + \cdots
x_{n}^2$.
We define the matrix
$$J_0 =  \left (
\begin{smallmatrix}
& &1  \\
& \adots &  \\\
1 & &
\end{smallmatrix} \right ) \;\;. $$
Let $\Lambda_0$ be the $(n+1) \times
(n+1)$ matrix $\textnormal{diag}(-1, 1, \cdots, 1)$ in the basis $(e_i)_{0 \leq i
  \leq n}$. Let $\Lambda'_0 = J_0 \cdot \Lambda_0 \cdot J_0 =
\textnormal{diag}(1, \cdots 1, -1)$. As a complex Hermitian form, $H$
is of signature $(2, 2n)$ with matrix
$ \Lambda = \textnormal{diag}(-1, 1 , \cdots, 1, -1)= \textnormal{diag}
(\Lambda_0, \Lambda'_0 ) $
in the basis $(f_i)_{0 \leq i \leq 2n+1}$.
Thus
$$ \Sp(n,1) = \{ X \in \GL(2n+2, \CC) \; / \; X \cdot
\left ( \begin{smallmatrix}
& J_1 \\
-J_1^t &
\end{smallmatrix} \right )
=
\left ( \begin{smallmatrix}
& J_1 \\
-J_1^t &
\end{smallmatrix} \right ) \cdot \overline{X}
\;\;  \textnormal{and} \;\; X^* \cdot \Lambda \cdot X
= \Lambda\;\;\}
 \;\;,$$
where $X^*$ denotes the complex trans-conjugate of $X$.
\begin{defi}
We denote by $J=
\left ( \begin{smallmatrix}
& J_0 \\
-J_0 &
\end{smallmatrix} \right )$ the product
$\Lambda \cdot \left ( \begin{smallmatrix}
& -J_1 \\
J_1 &
\end{smallmatrix} \right )
$.
\end{defi}
\noi
The complex symplectic form $\Omega$ on $\VC \oplus j \VC$ has matrix
$J$ in the basis $(f_i)_{0 \leq i \leq 2n+1}$. One can rewrite
$$ \Sp(n,1) = \{ X \in \GL(2n+2, \CC) \; / \;
X^* \cdot \Lambda \cdot X = \Lambda
\; \;  \textnormal{and} \;\; X^t \cdot J \cdot X
= J \;\;\}
 \;\;.$$
We thus recover the isomorphism $\Sp(n,1) = \U(2n, 2) \cap \Sp(2n+2,
\CC)$, where $\Sp(2n+2, \CC) = \Sp(\VC \oplus j \VC, \Omega)$.

\label{lie}
>From the previous descriptions we obtain~:
\begin{equation*}
\begin{split}
\LieSp(2n+2, \CC) &= \{
\left (
\begin{smallmatrix}
A & B \\
C & -J_0 \cdot A^t \cdot J_0
\end{smallmatrix} \right ) \; / \;J_0 \cdot C = C^t \cdot J_0 \;
\textnormal{and} \; J_0\cdot  B = B^t\cdot J_0 \} \;\;. \\
\LieU(2n,2) &= \{ \left (
\begin{smallmatrix}
A & B \\
- \Lambda_0 \cdot B^* \cdot \Lambda'_0 & D
\end{smallmatrix} \right )
\; / \; A^* \cdot \Lambda_0 + \Lambda_0 \cdot A =D^* \cdot \Lambda'_0 + \Lambda'_0
\cdot D=0  \}\;\;. \\
\LieSp(n,1)& = \{ \left (
\begin{smallmatrix}
A & B \\
- \Lambda_0 \cdot B^* \cdot \Lambda'_0 & -J_0 \cdot A^t \cdot J_0
\end{smallmatrix} \right ) \; /\; A^* \cdot \Lambda_0 + \Lambda_0
\cdot A = 0 = J_0 \cdot B - B^t \cdot J_0 \}\;\;.
\end{split}
\end{equation*}

Notice that the canonical embedding $ j_* :\su(n,1) \hookrightarrow \LieSp(n,1)$
factorizes through $\LieU(n,1)$. The embedding $\LieU(n,1)=
\{ A \; /\; A^* \cdot \Lambda_0 + \Lambda_0 \cdot A = 0 \}
 \lo \LieSp(n,1)$ is the morphism associating to $A \in \LieU(n,1)$
 the element $
\left ( \begin{smallmatrix}
A & 0 \\
0 & -J_0 \cdot A^t \cdot J_0
\end{smallmatrix} \right ) \in \LieSp(n,1)
$.

\subsubsection{Hodge filtration}
The action of the subgroup $S^1$ of $\Sp(n,1)$ with complexified Lie
algebra $\CC \cdot v$ (where
$v= \left (
\begin{smallmatrix}
-1 & 0 & 0 &0 \\
0& 0 &0 &0 \\
0 & 0 &0 &0 \\
0 & 0 & 0 & 1
\end{smallmatrix} \right ) \in \mathfrak{sp}(2n+2, \CC)$) on the Lie algebra $\LieSp(2n+2, \CC)$ defines a
filtration $F^\bullet \LieSp(2n+2, \CC)$ on $\LieSp(2n+2, \CC)$ with
graded pieces (we only indicate the non-positive ones)~:

\begin{align*}
\Gr^0 \LieSp(2n+2, \CC) &=\{
\left (
\begin{smallmatrix}
a &0 &0 &0 \\
0       &X &R &0\\
0       &S &-J_0 \cdot X^t \cdot J_0 &0 \\
0       &0 & 0 &-a
\end{smallmatrix} \right ) \; /
J_0\cdot  R = R ^t\cdot J_0 \; \text{and} \;J_0\cdot  S = S ^t\cdot
J_0 \} \;\;.  \\
\Gr^{-1}\LieSp(2n+2, \CC) &= \{
\left (
\begin{smallmatrix}
0        &u &v &0 \\
0       &0 &0 &v \cdot J_0\\
0       &0 &0 &u\cdot J_0 \\
0       &0 & 0 &0
\end{smallmatrix} \right ) u,v \in \CC^n \} \;\;. \\
\Gr^{-2} \LieSp(2n+2, \CC) &= \{
\left (
\begin{smallmatrix}
0 &0 &0 &z \\
0 &0 &0 &0\\
0 &0 &0 &0 \\
0 &0 & 0 &0
\end{smallmatrix} \right ) \; / z \in \CC\} \;\;.
\end{align*}

This Hodge decomposition restricts to a Hodge decomposition of the
complexified Lie algebra $\gl(n+1, \CC) = \su(n,1) \otimes_\RR \CC$.

\subsubsection{Automorphic bundles}
Let $K$ denote the maximal compact subgroup $S(U(n) \times U(1))$ of
$SU(n,1)$ and $K_\CC \simeq GL(n, \CC)$ its complexification. This is
a Levi subgroup of the parabolic subgroup $Q \cap \SU(n,1)(\CC)$,
where $Q$ denotes the parabolic subgroup of $\Sp(2n+2, \CC)$ with Lie
algebra $F^0\LieSp(2n+2, \CC)$. The
natural inclusion
$$ \HypC^n = SU(n,1)/K \hookrightarrow \SU(n,1)(\CC)/(Q \cap
\SU(n,1)(\CC)) \simeq \mathbf{P}^n \CC$$ is the natural open embedding of
the period domain $\HypC^n$ into its dual.

\begin{defi}
Given a $K_\CC$-module $\mathfrak{m}$, we denote by $\F(\mathfrak{m})$ the
holomorphic automorphic vector bundle $\Gamma \backslash (\SU(n,1)(\CC)\times_{Q
  \cap \SU(n,1)(\CC)} \mathfrak{m})_{|\HypC^n}$ with fiber $\mathfrak{m}$ on $M$.
\end{defi}

\begin{defi}
We denote by $L_M$ the automorphic line bundle $n+1$-th root of $K_M^{-1}$, where $K_M$ denotes the canonical line bundle on $M$.
\end{defi}

\subsubsection{Non-Abelian Hodge theory}

Let $(\Gr \mathcal{P}_f, \theta_f)$ be the system of Hodge bundles associated to
the variation of Hodge structures $(\rho, f)$. Thus $$\Gr
\mathcal{P}_f = \F(\Gr^\bullet \LieSp(2n+2, \CC))\;\;.$$
As proven in \cite{klingler} the
tangent space $T$ to the space of $\Sp(n,1)$-variations of Hodge structures
(equivalently~: to the subspace of systems of Hodge bundles in the
space of semistable $\textnormal{Sp}(2n+2, \CC)$-Higgs bundles) on
$M$, modulo the trivial deformations in the centralizer $U(n,1)$ of $SU(n,1)$,
identifies with the hypercohomology of complex of coherent sheaves~:
\begin{equation*}
\HH^1(M, \F(\frac{\Gr^0 \LieSp(2n+2, \CC)}{\Gr^0 \gl(n+1, \CC)})
\stackrel{\theta_f}{\lo} \F(\frac{\Gr^{-1} \LieSp(2n+2, \CC)}{\Gr^{-1}
  \gl(n+1, \CC)})
  \otimes \Omega^1_M  \stackrel{\theta_f}{\lo} \F(\frac{\Gr^{-2}
    \LieSp(2n+2, \CC)}{\Gr^{-2} \gl(n+1, \CC)})
  \otimes \Omega^2_M )\;\;.
\end{equation*}
One easily computes~:
\begin{align*}
\F(\frac{\Gr^0 \LieSp(2n+2, \CC)}{\Gr^0 \gl(n+1, \CC)})&=
\left (
\begin{smallmatrix}
0 &0 &0 &0 \\
0       &0 &S^2 \Omega^1_M \otimes L_M^{2} &0\\
0       &S^2 TM \otimes L_M^{-2} &0 &0 \\
0       &0 & 0 &0
\end{smallmatrix} \right ) \\
\F(\frac{\Gr^{-1} \LieSp(2n+2, \CC)}{\Gr^{-1}
  \gl(n+1, \CC)}) &=
\left (
\begin{smallmatrix}
0        &0 &v \in \Omega^1_M \otimes L_M^{2} &0 \\
0       &0 &0 &v \cdot J_0\\
0       &0 &0 &0 \\
0       &0 & 0 &0
\end{smallmatrix} \right )  \\
\F(\frac{\Gr^{-2}
    \LieSp(2n+2, \CC)}{\Gr^{-2} \gl(n+1, \CC)}) &=
\left (
\begin{smallmatrix}
0 &0 &0 &L_M^{2} \\
0 &0 &0 &0\\
0 &0 &0 &0 \\
0 &0 & 0 &0
\end{smallmatrix}\right ) \;\;.
\end{align*}

Thus
$$
T=\HH^1(M,
\left (
\begin{smallmatrix}
0 &S^2 \Omega^1_M \otimes L_M^{2} \\
S^2 T_M \otimes L_M^{-2} & 0
\end{smallmatrix} \right )
\stackrel{\left ( \begin{smallmatrix} 0 \\ 1 \end{smallmatrix} \right )}{\lo}
\left(  \begin{smallmatrix} (\Omega^1_M \otimes L_M^{2})\otimes
    \Omega^1_M \\0
  \end{smallmatrix} \right )
\stackrel{\left ( \begin{smallmatrix} 0 \\ 1 \end{smallmatrix} \right )}{\lo} L_M^{2} \otimes \Omega^2_M )\;\;.$$
Thus~:
\begin{align*}
\begin{split}
T &=\HH^1(M, \left (
\left ( \begin{smallmatrix} S^2T_M \otimes L_M^{-2} \\ S^2\Omega^1_M \otimes L_M^{2} \end{smallmatrix}\right )
\stackrel{\left ( \begin{smallmatrix} 0 &0 \\ 1 & 0 \end{smallmatrix} \right )}{\lo}
\left ( \begin{smallmatrix} S^2\Omega^1_M \otimes L_M^{2} \\ \Omega^2_M \otimes L_M^{-2} \end{smallmatrix}\right )
\stackrel{\left ( \begin{smallmatrix} 0 &1 \end{smallmatrix} \right )}{\lo}
\Omega^2_M \otimes L_M^{2}
\right ) \\
&= H^1(M,S^2T_M \otimes L_M^{-2})\;\;.
\end{split}
\end{align*}
Notice that $H^1(M,S^2T_M \otimes L_M^{-2})$ is
conjugate to $H^0(M, S^3 \Omega^1_M\otimes L_M^2)$ via the natural
pairing~:
$$ H^0(M, S^3 \Omega^1_M \otimes L_M^{2}) \otimes H^1(M,S^2T_M \otimes
L_M^{-2}) \lo H^1(M, \Omega^1_M)= H^{1,1}(M, \CC)
\stackrel{\cdot \wedge \omega_M^{n-1}}{\lo} \CC\;\;.$$
This finishes the proof of proposition~\ref{pi'}.

\sspace
\noi

\noi
Inkang Kim\\
School of Mathematics\\
     KIAS, Hoegiro 87, Dongdaemun-gu\\
     Seoul, 130-722, Korea\\
     \texttt{inkang\char`\@ kias.re.kr}

\sspace
\noi Bruno Klingler\\
Institut de Math\'ematiques de Jussieu, Paris, France \\
and Institute for Advanced Study, Princeton NJ, USA \\
\texttt{klingler\char`\@ math.jussieu.fr}

\sspace
\noi Pierre Pansu\\ Laboratoire de Math{\'e}matiques
     d'Orsay\\
     UMR 8628 du CNRS\\
 Universit{\'e} Paris-Sud\\
 91405 Orsay C\'edex, France\\
  \texttt{pierre.pansu\char`\@ math.u-psud.fr}

\end{document}